\documentclass[a4paper,12pt]{article}
\usepackage{a4wide}
\usepackage{amsmath}
\usepackage{amssymb}
\usepackage{amsthm}
\usepackage{latexsym}
\usepackage{graphicx}
\usepackage[english]{babel}
\usepackage{makeidx}
              
\newtheorem{dfn} [subsection]{Definition}
\newtheorem{obs} [subsection]{Remark}
\newtheorem{exm} [subsection]{Example}

\newtheorem{prop}[subsection]{Proposition}

\newtheorem{teor}[subsection]{Theorem}
\newtheorem{lema}[subsection]{Lemma}
\newtheorem{cor} [subsection]{Corollary}

\newcommand{\de}{\mathbf{d}}
\newcommand{\me}{\mathbf{m}}

\newcommand{\Dei}{$\mathcal D$-fixed ideal }
\newcommand{\Deis}{$\mathcal D$-fixed ideals }
\newcommand{\Deisp}{$\mathcal D$-fixed ideals.}
\newcommand{\Deisv}{$\mathcal D$-fixed ideals, }
\newcommand{\Di}{$\mathcal D$-fixed ideal }

\newcommand{\di}{$\de$-fixed ideal }
\newcommand{\dis}{$\de$-fixed ideals }

\begin{document}

\selectlanguage{english}
\frenchspacing

\large
\begin{center}
\textbf{Monomial ideals with linear upper bound regularity.}

Mircea Cimpoea\c s
\end{center}

\normalsize

\footnotetext[1]{This paper was supported by the CEEX Program of the Romanian
Ministry of Education and Research, Contract CEX05-D11-11/2005 and by 
 the Higher Education Commission of Pakistan.}

\begin{abstract}
In this paper, we extend a result of Eisenbud-Reeves-Totaro in the frame of ideals of Borel type. 
As a consequence, we obtain a linear upper bound for the regularity of a new class of ideals, called \Deisp

\vspace{5 pt} \noindent \textbf{Keywords:} Borel type ideals, $p$-Borel ideals, Mumford-Castelnuovo regularity.

\vspace{5 pt} \noindent \textbf{2000 Mathematics Subject
Classification:}Primary: 13P10, Secondary: 13E10.
\end{abstract}

\begin{center}
\textbf{Introduction.}
\end{center}

Let $K$ be an infinite field, and let $S=K[x_1,...,x_n],n\geq 2$ the polynomial ring over $K$.
Bayer and Stillman \cite{BS} note that a Borel fixed ideals $I\subset S$ satisfy the following property:
\[(*)\;\;\;\;(I:x_j^\infty)=(I:(x_1,\ldots,x_j)^\infty)\;for\; all\; j=1,\ldots,n.\] Herzog, Popescu and Vladoiu \cite{hpv} say that a monomial ideal $I$ is of \emph{Borel type}, if it satisfy $(*)$. We mention that this concept appears also in \cite{CS} as the so called weakly stable ideal. Herzog, Popescu and Vladoiu proved in \cite{hpv} that $I$ is of Borel type, if and only if for any monomial $u\in I$ and for any $1\leq j<i \leq n$ and $q>0$ with $x_i^{q}|u$, there exists an integer $t>0$ such that $x_j^{t}u/x_i^{q}\in I$. In the first section of our paper, we prove that if $I$ is an ideal of Borel type, then $I_{\geq e} = $ the ideal generated by the monomials of degree $\geq e$ from $I$, is stable whenever $e\geq reg(I)$ (Theorem $1.5$). This allow us to give a generalization of a result of Eisenbud-Reeves-Totaro (Corollary $1.8$). 

In the second section, we introduce a new class of monomial ideals, called \Deisv which are sums of various \dis, introduced in \cite{mir} (see also Section $2$). A.Imran and A.Sarfraz show in \cite{saf} that an ideal of Borel type  $I\subset S$ has the regularity upper bounded by $n \cdot (deg(I)-1)+1$, where $deg(I)$ is the maximum degree of a minimal generator of $I$. In particular, this holds for \Deis, since any \Dei is a Borel type ideal (Proposition $2.3$). 

A formula for the regularity of a principal \di, i.e a the smallest \di which contain a monomial $u\in S$, was given in \cite[Theorem 3.1]{mir}. This result together with Corollary $1.6$ allow us to give a new proof to the fact that a \Dei $I$ has the regularity upper bounded by $n \cdot (deg(I)-1)+1$ (Corollary $2.5$). Also, Example $2.7$ show that the family of all \dis is strictly included in the family of \Deis, which is, also, strictly included in the family of
ideals of Borel type.

The author is grateful to his adviser Dorin Popescu for his encouragement and valuable suggestions.
My thanks goes also to the School of Mathematical Sciences, GC University, Lahore, Pakistan for supporting 
and facilitating this research.

\newpage
\section{A stable property of Borel type ideals.}

\begin{obs}
{\em If $I,J$ are two ideals of Borel type then $I+J$, $I\cap J$ and $I\cdot J$ are of Borel type. Also, a quotient ideal of an ideal of Borel type by an monomial ideal is of Borel type.}
\end{obs}

\begin{obs}{\em
For any monomial $u\in S$ we denote $m(u)=max\{i:\; x_i|u\}$. For any monomial ideal $I\subset S$, we denote
$m(I)=max\{m(u):\;u\in G(I)\}$, where $G(I)$ is the set of the minimal generators of $I$. Also, if $M$ is a graded $S$-module of finite length, we denote $s(M)=max\{t:\;M_{t}\neq 0\}$.

Let $I\subset S$ a Borel type ideal. In \cite{hpv}, it was defined a chains of ideals as follows. We let $I_0=I$. Suppose $I_{\ell}$ is already defined. If $I_{\ell}=S$ then the chain ends. Otherwise, we let $n_{\ell}=m(I_{\ell})$ and set $I_{\ell+1}=(I_{\ell}:x_{n_{\ell}}^{\infty})$. Notice that $r\leq n$, since $n_{\ell}>n_{\ell+1}$ for all $0\leq \ell<r$. The chain $I = I_{0} \subset I_{1}\subset \cdots \subset I_r = S$ is called the \emph{sequential chain} of $I$.
\cite[Corollary 2.5]{hpv} states: 
\[(1)\; I_{\ell+1}/I_{\ell} \cong (J_{\ell}^{sat}/J_{\ell})[x_{n_{\ell}+1},\ldots,x_n],\]
for all $0\leq \ell<r$, where $J_{\ell}\subset S_{\ell} = K[x_1,\ldots,x_{n_{\ell}}]$ is the ideal generated by $G(I_{\ell})$. Also, using $(1)$, \cite[Corollary 2.5]{hpv} gives a formula for the Mumford-Castelnuovo regularity of $I$, more precisely, \[(2)\; reg(I) = max\{s(J_{0}^{sat}/J_0),s(J_{1}^{sat}/J_1),\cdots,s(J_{r-1}^{sat}/J_{r-1})\} + 1.\]
}
\end{obs}

\begin{lema}
Let $I\subset S$ be an ideal and $I'=IS'$ the extension of $I$ in $S'=S[x_{n+1}]$. If $e\geq deg(I)$, then
$I_{\geq e}$ is stable if and only if $I'_{\geq e}$ is stable.
\end{lema}

\begin{proof}
Let $u\in I'_{\geq e}$ be a monomial. Then $u=x_{n+1}^{k}\cdot v$ for some $v\in I$. If $k>0$ then $m(u)=n+1$ and therefore, for any $i<n+1$, $x_i\cdot u/x_{n+1} = x_{n+1}^{k-1}\cdot x_i\cdot v \in I'_{\geq e}$. If $k=0$ then
$m(u)\leq n$ and since $I_{\geq e}$ is stable, it follows $x_i \cdot u/ x_{m(u)} \in I'_{\geq e}$. Thus $I'_{\geq e}$
is stable. For the converse, simply notice that $G(I_{\geq e})\subset G(I'_{\geq e})$ and since is enough to
check the stable property only for the minimal generators, we are done.
\end{proof}

\begin{lema}
If $I\subset S$ is an artinian monomial ideal and $e\geq reg(I)$ then $I_{\geq e}$ is stable.
\end{lema}

\begin{proof}
Note that $reg(I) = s(S/I) + 1$, therefore, if $e\geq reg(I)$ then $I_{\geq e} = S_{\geq e}$, thus $I_{\geq e}$ is stable.
\end{proof}

\begin{teor}
Let $I\subset S$ be a Borel type ideal and let $e\geq reg(I)$ be an integer. Then $I_{\geq e}$ is stable.
\end{teor}

\begin{proof}
We use induction on $r\geq 1$, where $r$ is the length of the sequential chain of $I$. If $r=1$, i.e. $I$ is an
artinian ideal, we are done from the previous lemma. 

Suppose now $r>1$ and let $I = I_0 \subset I_1 \subset \cdots \subset I_r = S$ be the sequential chain of $I$. Let $n_{\ell}=m(I_{\ell})$ for  $0\leq \ell \leq r$. Let $J_{\ell}\subset S_{\ell} = K[x_{1},\ldots,x_{n_{\ell}}]$ be the ideal generated by $G(I_{\ell})$. Using the induction hypothesis, we may assume $(I_1)_{\geq e}$ stable for $e\geq reg(I_1)$. On the other hand, from $(2)$ it follows that $reg(I_1)\leq reg(I)$, thus 
$(I_1)_{\geq e}$ is stable for $e\geq reg(I)$.

Since $J_0^{sat} = I_1 \cap S_{n_0}$, using Lemma $1.4$ it follows that $(J_0^{sat})_{\geq e}$ is stable.
Let $e\geq reg(I)$. Since $reg(I) \geq s(J_0^{sat}/J_0) + 1$ it follows
that $(J_0)_{\geq e} = (J_0^{sat})_{\geq e}$ is stable. Since $I_0 = J_0S$, using again Lemma $1.4$, we get
$I_{\geq e}$ stable for $e\geq reg(I)$, as required.
\end{proof}

\begin{prop}\cite[Proposition 12]{ert}
Let $I$ be a monomial ideal with $deg(I)=d$ and let $e\geq d$ such that $I_{\geq e}$ is stable. Then $reg(I)\leq e$.
\end{prop}

Theorem $1.6$ and Proposition $1.7$ yield the following:

\begin{cor}
If $I$ is a Borel type ideal, then $reg(I)=min\{e:\; I_{\geq e}\;is\; stable \}$.
\end{cor}

\begin{cor}
If $I$ is a monomial ideal with $Ass(S/I)$ totally ordered by inclusion, then $reg(I)=min\{e:\; I_{\geq e}\;is\; stable \}$.
\end{cor}

\begin{proof}
By renumbering the variables, we can assume that each $p\in Ass(S/I)$ has the form 
$p=(x_1,x_2,...,x_r)$ for some $1\leq r\leq n$. Therefore, by \cite[Theorem 2.2]{saf}, we may assume
that $I$ is an ideal of Borel type, and then we apply Corollary $1.8$.
\end{proof}

\begin{exm}{\em
Let $I = (x_1^7,\; x_1^5x_2,\; x_1^2x_2^4,\; x_1x_2^6,\; x_1^5x_3^2,\; x_1x_2^4x_3^2)\subset K[x_1,x_2,x_3,x_4]$.
We construct the sequential chain of $I$. We have $I_0 = I$ and $n_0=m(I_0)=3$, therefore $J_0 = I_0 \cap K[x_1,x_2,x_3]$. Let $I_1 = (I_0 : x_3^{\infty}) = (x_1^5,x_1x_2^4)$. We have $n_1 = m(I_1) = 2$, therefore
$J_1 = I_1 \cap K[x_1,x_2]$. Let $I_2 = (I_1 : x_2^{\infty}) = (x_1)$. We have $n_2 = m(I_2) = 1$, therefore
$J_2 = I_2 \cap K[x_1]$. One can easily compute, $s(J_0^{sat}/J_0) = 7$, $s(J_1^{sat}/J_1)=7$ and $s(J_2^{sat}/J_2)=1$.
Using $(2)$, we get:
\[ reg(I) = max\{s(J_0^{sat}/J_0), s(J_1^{sat}/J_1), s(J_2^{sat}/J_2)\} + 1 = 8. \]
We will exemplify the proof of $1.5$ for $I$. Let $e\geq 8$ be an integer. Since $(J_2)_{\geq e}$ is obviously stable, it follows from $1.3$ that $(I_2)_{\geq e} = (J_2S)_{\geq e}$ is also stable. Note that $I_2 = J_1^{sat}S$ and moreover, $I_2$ and $J_1^{sat}$ have the same minimal set of generators. Therefore, by $1.3$, it follows that $(J_1^{sat})_{\geq e}$ is stable. Since $e\geq reg(I)> s(J_1^{sat}/J_1) = 7$, it follows that $(J_1)_{\geq e}$ is stable.
On the other hand, $I_1=J_1S$, therefore $(I_1)_{\geq e}$ is stable. Since $J_0^{sat}S = I_1$ we get, from $1.3$,
$(J_0^{sat})_{\geq}$ stable. Since $e\geq reg(I)> s(J_0^{sat}/J_0) = 7$, it follows that $(J_0)_{\geq e}$ is stable.
Finally, since $I = I_0 = J_0S$, we obtain $I_{\geq e}$ stable, as required.
}\end{exm}

\section{\Deis}

A \emph{$\de$-sequence} is a strictly increasing sequence of integers, $\de: 1=d_{0}|d_{1}|\cdots|d_{s}$. For any $a\in \mathbb N$ there exists an unique sequence of positive integers $a_{0},a_{1},\ldots,a_{s}$ such that $a= \sum_{t=0}^{s}a_{t}d_{t}$ and $0 \leq a_{t} < \frac{d_{t+1}}{d_{t}}$, for any $0\leq t<s$.
The converse is also true, in the sense that if $1=d_{0}<d_{1}<\cdots<d_{s}$ is a sequence of positive integers such that for any $a \in \mathbb N$ there exists an unique sequence of positive integers $a_{0},a_{1},\ldots,a_{s}$ as before, then $1=d_{0}|d_{1}|\cdots|d_{s}$ (see \cite[Lemma 1.1]{mir}). The decomposition $a=\sum_{t=0}^{s}a_{t}d_{t}$ is called the \emph{$\de$-decomposition} of $a$. In particular, if $d_{t}=p^{t}$ we get the $p$-adic decomposition of $a$. 

Let $a,b\in\mathbb N$ and consider the decompositions $a=\sum_{t=0}^{s}a_{t}d_{t}$ and $b=\sum_{j=0}^{s}b_{t}d_{t}$. 
We write $a\leq_{\de}b$ if $a_{t}\leq b_{t}$ for any $0\leq t\leq s$. This allowed us to define the following class of monomial ideals:
We say that a monomial ideal $I\subset S$ is \emph{$\de$-fixed}, if for any monomial $u\in I$ and for any indices $1\leq j<i \leq n$, if $t\leq_{\de} \nu_{i}(u)$ (where $\nu_{i}(u)$ is the exponent of the variable $x_{i}$ in $u$) then $u \cdot x_{j}^{t}/x_{i}^{t}\in I$
(see \cite[Definition 1.4]{mir}).

Let $2\leq i_{1}<i_{2}<\cdots<i_{r}=n$ and let $\alpha_{1},\ldots,\alpha_{r}$ some positive integers.
Let $u=\prod_{i=1}^{r}x_{i_{q}}^{\alpha_{q}} \in S$. Let $I = <u>_{\de}$ the principal
\di generated by $u$. From \cite[Proposition 1.8]{mir} it follows that 
$ I = \prod_{r=1}^{q}\prod_{j=0}^{s}(\me_{q}^{[d_{j}]})^{\alpha_{qj}}$, 
where $\alpha_{q}=\sum_{j=0}^{s}\alpha_{qj}d_{j}$. Let $s_{q}=max\{j|\alpha_{qj}\neq 0\}$ and $d_{qt}= \sum_{e=1}^{q}\sum_{j\geq t}^{s_{e}}\alpha_{ej}d_{j}$. Let $D_{q} = d_{qs_{q}} + (i_{q}-1)(d_{s_{q}}-1)$, for $1\leq q \leq r$. With this notations we have, the following generalization of so called Pardue's formula (see \cite{ah}, \cite{hp} and \cite{hpv}):

\begin{teor}(\cite[Theorem 3.1]{mir})
$reg(I) = max_{1\leq q \leq r} D_{q}$. In particular, if $I = <x_{n}^{\alpha}>_{\de}$ and $\alpha=\sum_{t=0}^{s} \alpha_{t}d_{t}$ with $\alpha_{s}\neq 0$ then $reg(I) = \alpha_{s}d_{s} + (n-1)(d_{s}-1)$.
\end{teor}

\begin{dfn}
We say that a monomial ideal $I\subset S$ is a \Dei, if $I$ is a sum of $\de$-fixed ideals, 
for various $\de$-sequences.
\end{dfn}

Since any $\de$-fixed ideal is a Borel type ideal (see \cite[Proposition 1.10]{mir}), we easily get 
from Remark $1.1$ the following:

\begin{prop}
If $I$ is a \Dei then $I$ is a Borel type ideal.
\end{prop}

Therefore, from Corollary $1.8$ we get the next:

\begin{cor}
If $I$ is a \Di then $reg(I)=min\{e:\; I_{\geq e}\;is\; stable \}$.
\end{cor}

We mention that this result was first obtained as a consequence of the proof of Pardue's formula by Herzog-Popescu \cite{hp} in the special case of a principal $p$-Borel ideal.

\begin{cor}
Let $I$ be a \Dei. Suppose $I = I_{1} + \cdots + I_{m}$ where $I_{i}$ are principal $\de^{(i)}$-fixed ideals for several
$\de$-sequences, $\de^{(i)}$ where $i=1,\ldots,m$. Then
\[ reg(I) \leq max\{reg(I_{1}),\ldots,reg(I_{m})\}.\]
\end{cor}

\begin{proof}
From Corollary $2.4$ it follows that $(I_{i})_{\geq reg(I_{i})}$ is stable for $i=1,\ldots,r$. So, if $e = max\{reg(I_{1}),\ldots,reg(I_{r})\}$ it follows that $(I_{i})_{\geq e}$ are stables for all $i$. Therefore
$I_{\geq e}=(\sum_{i=1}^{r}I_{i})_{\geq e} = \sum_{i=1}^{r}(I_{i})_{\geq e}$ is stable, as a sum of stable ideals.
Thus, from $2.4$ we get $reg(I)\leq e$ as required.
\end{proof}

It was already proved by A. Imran and A. Sarfraz in \cite{saf} that for a Borel type ideal $I$, \linebreak
$reg(I)\leq n \cdot (deg(I)-1)+1$. Using Theorem $2.1$ and the previous corollary, we can give another proof of this result in the particular case of \Deisp

\begin{cor}
Let $I$ be a $\mathcal D$-fixed ideal. Then $reg(I)\leq n \cdot (deg(I)-1)+1$. In particular, if $e\geq n \cdot (deg(I)-1)+1$ then $I_{\geq e}$ is stable.
\end{cor}

\begin{proof}
Suppose $I = I_{1} + \cdots + I_{m}$, where $I_{i}$ are principal $\de^{(i)}$-fixed ideals for several \linebreak
$\de$-sequences, $\de^{(i)}$ with $i=1,\ldots,m$. From Theorem $2.1$ one can easily see that \linebreak
$reg(I_{i}) \leq n\cdot (deg(I)-1) +1$ for all $1\leq i\leq m$, thus, using $2.5$, we are done.
\end{proof}


\begin{exm}{\em
	(1)  Let $d_1,d_2,\cdots,d_n>0$ be some integers. Then $I=(x_1^{d_1},x_2^{d_2},\ldots,x_n^{d_n})$ is a \Di if
	      and only if $d_1\leq d_2 \leq \cdots \leq d_n$. Indeed, if there exists a $j<i$ such that $d_j>d_i$, since
	      $I$ is a \Di and $x_i^{d_i}\in I$ it follows $x_j^{d_i}\in I$ which is a contradiction. Indeed, it is enough
	      to note that $x_i^{d_i} \in I'\subset I$, where $I'$ is a \di and then we apply the definition of a \di.
	      For the converse, we consider $\mathbf{d}^{(i)}$
	      to be the $\de$-sequence $1|d_i$. We have $<x_i^{d_i}>_{\mathbf{d}^{(i)}} = (x_1^{d_i},\ldots,x_i^{d_i})$ thus
	      $I=<x_1^{d_1}>_{\mathbf{d}^{(1)}} + \cdots + <x_n^{d_n}>_{\mathbf{d}^{(n)}}$. Since any artinian monomial ideal
	      is a Borel type ideal, this example shows that a Borel type ideal is not in general a \Di. Modulo renumbering 
	      the variables, any monomial complete intersection is a \Dei. But not any artinian ideal can be obtained from
	      a \Dei in such way, as example $(x_1^{3},x_1x_2, x_3^2)$ shows. 
	      
	(2) In general, if $0<d_1\leq d_2 \leq \cdots \leq d_n$ are some integers,
	      $I=(x_1^{d_1},x_2^{d_2},\ldots,x_n^{d_n})$ is not a \di. Indeed, let
	      $I=(x_1^{2},x_2^{3},x_3^{4})$ and let $\de:1=d_0|d_1|d_2|\cdots$ be a $\de$-sequence. If $d_1>4$ then, since
	      $x_3^{4}\in I$ and $4=4\cdot d_0$ it follows that $(x_1,x_2,x_3)^{4}\subset I$, a contradiction. If $d_1=4$,
	      then, since $x_2^{3}\in I$ it follows that $(x_1,x_2)^{3}\subset I$, a contradiction. If $d_1=3$, 
	      since $x_3^{4}\in I$ and $4=1+3$, it follows that $x_3^{3}x_2\in I$, a contradiction. If $d_1=2$,
	      since $x_2^{3}\in I$ and $3=1+2$, it follows that $x_1x_2^{2}\in I$, a contradiction. Therefore, $I$
	      is not a \di for any $\de$-sequence.
}\end{exm}

\vspace{2mm} \noindent {\footnotesize
\begin{minipage}[b]{15cm}
 Mircea Cimpoea\c s, Institute of Mathematics of the Romanian Academy, Bucharest, Romania\\
 E-mail: mircea.cimpoeas@imar.ro
\end{minipage}}
\end{document}